\documentclass[amstex,12pt, amssymb]{article}

\usepackage{mathtext}
\usepackage[cp1251]{inputenc}
\usepackage[T2A]{fontenc}
\usepackage[dvips]{graphicx}
\usepackage{amsmath}
\usepackage{amssymb}
\usepackage{amsxtra}
\usepackage{latexsym}
\usepackage{ifthen}

\newcounter{lemma}[section]

\newcounter{corollary}[section]

\newcounter{remark}[section]

\newcounter{theorem}[section]

\newcounter{proposition}[section]

\numberwithin{equation}{section}

\def\Xint#1{\mathchoice
   {\XXint\displaystyle\textstyle{#1}}%
   {\XXint\textstyle\scriptstyle{#1}}%
   {\XXint\scriptstyle\scriptscriptstyle{#1}}%
   {\XXint\scriptscriptstyle\scriptscriptstyle{#1}}%
   \!\int}
\def\XXint#1#2#3{{\setbox0=\hbox{$#1{#2#3}{\int}$}
     \vcenter{\hbox{$#2#3$}}\kern-.5\wd0}}
\def\dashint{\Xint-}

\def\Rn{{{\Bbb R}^n}}

\def\cc{\setcounter{equation}{0}
\setcounter{figure}{0}\setcounter{table}{0}}

\begin{document}

\markboth{\centerline{DENIS KOVTONYUK AND VLADIMIR RYAZANOV}}
{\centerline{ON BOUNDARY BEHAVIOR OF SPATIAL MAPPINGS}}

\author{{DENIS KOVTONYUK AND VLADIMIR RYAZANOV}}

\title{{\bf ON BOUNDARY BEHAVIOR\\ OF SPATIAL MAPPINGS}}

\maketitle

\large \begin{abstract} We show that homeomorphisms $f$ in ${\Bbb
R}^n$, $n\geqslant3$, of finite distortion in the Orlicz--Sobolev
classes $W^{1,\varphi}_{\rm loc}$ with a condition on $\varphi$ of
the Calderon type and, in particular, in the Sobolev classes
$W^{1,p}_{\rm loc}$ for $p>n-1$ are the so-called lower
$Q$-homeomorphisms, $Q(x)=K^{\frac{1}{n-1}}_I(x,f)$, where
$K_I(x,f)$ is its inner dilatation. The statement is valid also for
all finitely bi-Lipschitz mappings that a far--reaching extension of
the well-known classes of isometric and quasiisometric mappings.
This makes pos\-sib\-le to apply our theory of the boundary behavior
of the lower $Q$-homeomorphisms to all given classes.
\end{abstract}


\large \cc
\section{Introduction}

The problem of the boundary behavior is one of the central topics of
the theory of quasiconformal mappings and their generalizations. At
present mappings with finite distortion are  studied, see many
references in the monographs \cite{GRSY} and \cite{MRSY}. In this
case, as it was earlier, the main geometric approach in the modern
mapping theory is the method of moduli, see, e.g., the monographs
\cite{GRSY}, \cite{MRSY}, \cite{Oht}, \cite{Ri}, \cite{Vs},
\cite{Va}  and \cite{Vu}.

\medskip

It is well--known that the concept of moduli with weights
essentially due to Andreian Cazacu, see, e.g.,
\cite{And1}--\cite{AC$_1$}, see also recent works
\cite{Cr1}--\cite{Cr3} of her learner. At the present paper we give
new modulus estimates for space mappings that essentially improve
the corresponding estimates first obtained in the paper
\cite{KRSS1}.

\medskip

Here we apply our theory of the so--called lower $Q$-homeomorphisms
first introduced in the paper \cite{KR$_1$}, see also the monograph
\cite{MRSY}, that was motivated by the ring definition of
quasiconformal mappings of Gehring, see \cite{Ge$_1$}. The theory of
lower $Q$-homeomorphisms has already found interesting applications
to the theory of the Beltrami equations in the plane and to the
theory of mappings of the classes of Sobolev and Orlich-Sobolev in
the space, see, e.g., \cite{KPR},
 \cite{KPRS}, \cite{KRSS1}, \cite{KRSS}, \cite{KSS}, \cite{MRSY} and
\cite{RSSY}.

\medskip

Following Orlicz, see, e.g., paper \cite{Or1}, see also monograph
\cite{Za}, given a convex increasing function $\varphi:{\Bbb
R}^+\to{\Bbb R}^+$, $\varphi(0)=0$, denote by $L^{\varphi}$ the
space of all functions $f:D\to{\Bbb R}$ such that
\begin{equation}\label{eqOS1.1}
\int\limits_{D}\varphi\left(\frac{|f(x)|}{\lambda}\right)\ dm(x)\ <\
\infty\end{equation} for some $\lambda>0$ where $dm(x)$ corresponds
to the Lebesgue measure in $D$. $L^{\varphi}$ is called the {\bf
Orlicz space}. In other words, $L^{\varphi}$ is the cone over the
class of all functions $g:D\to{\Bbb R}$ such that
\begin{equation}\label{eqOS1.2}
\int\limits_{D}\varphi\left(|g(x)|\right)\ dm(x)\ <\
\infty\end{equation} which is also called the {\bf Orlicz class},
see \cite{BO}.

\medskip

The {\bf Orlicz--Sobolev class} $W^{1,\varphi}(D)$ is the class of
all functions $f\in L^1(D)$ with the first distributional
derivatives whose gradient $\nabla f$ belongs to the Orlicz class in
$D$. $f\in W^{1,\varphi}_{\rm loc}(D)$ if $f\in W^{1,\varphi}(D_*)$
for every domain $D_*$ with a compact closure in $D$. Note that by
definition $W^{1,\varphi}_{\rm loc}\subseteq W^{1,1}_{\rm loc}$. As
usual, we write $f\in W^{1,p}_{\rm loc}$ if $\varphi(t)=t^p$,
$p\geqslant1$. Later on, we also write $f\in W^{1,\varphi}_{\rm
loc}$ for a locally integrable vector-function $f=(f_1,\ldots,f_m)$
of $n$ real variables $x_1,\ldots,x_n$ if $f_i\in W^{1,1}_{\rm loc}$
and
\begin{equation}\label{eqOS1.2a} \int\limits_{D}\varphi\left(|\nabla
f(x)|\right)\ dm(x)\ <\ \infty\end{equation} where $|\nabla
f(x)|=\sqrt{\sum\limits_{i,j}\left(\frac{\partial f_i}{\partial
x_j}\right)^2}$. Note that in this paper we use the notation
$W^{1,\varphi}_{\rm loc}$ for more general functions $\varphi$ than
in those classic Orlicz classes often giving up the conditions on
convexity and normalization of $\varphi$. Note also that the
Orlicz--Sobolev classes are intensively studied in various aspects
at the moment, see, e.g., \cite{KRSS} and references therein.

\medskip

In this connection, recall definitions which are relative to
Sobolev's classes. Given an open set $U$ in $\mathbb{R}^n$, $n\ge
2,$ $C_0^{ \infty }(U)$ denotes the collection of all functions
$\psi : U \to \mathbb{R}$ with compact support having continuous
partial derivatives of any order. Now, let $u$ and $v: U \to {\Bbb
R}$ be locally integrable functions. The function $v$ is called the
{\bf distributional (generalized) derivative} $u_{x_i}$ of $u$ in
the variable $x_i$, $i=1,2,\ldots , n$, $x=(x_1,x_2,\ldots , x_n)$,
if
\begin{equation}\label{eqSTR2.21}
\int\limits_{U} u\, \psi_{x_i} \,dm(x)=- \int\limits_{U} v\, \psi \
dm(x) \ \quad \forall \ \psi \in C_{0}^{\infty}(U)\ .
\end{equation}
$u\in W^{1,1}_{\rm loc}$ if $u_{x_i}$ exist for all $i=1,2,\ldots ,
n$. These concepts were introduced by Sobolev in ${\Bbb R}^n$,
$n\geqslant2$, see \cite{So}, and at present it is developed under
wider settings by many authors, see, e.g., \cite{Re$_1$} and further
references in \cite{KRSS}.

\medskip

Recall also the following topological notion. A domain
$D\subset{\Bbb R}^n$, $n\geqslant2$, is said to be {\bf locally
connected at a point} $x_0\in\partial D$ if, for every neighborhood
$U$ of the point $x_0$, there is a neighborhood $V\subseteq U$ of
$x_0$ such that $V\cap D$ is connected. Note that every Jordan
domain $D$ in ${\Bbb R}^n$ is locally connected at each point of
$\partial D$, see, e.g., \cite{Wi}, p. 66.

Following \cite{KR} and \cite{KR$_1$}, see also \cite{MRSY} and
\cite{RSal}, we say that $\partial D$ is {\bf weakly flat at a
point} $x_0\in\partial D$ if, for every neighborhood $U$ of the
point $x_0$ and every number $P>0$, there is a neighborhood
$V\subset U$ of $x_0$ such that
\begin{equation}\label{eq1.5KR} M(\Delta(E,F,D))\ \geqslant P\end{equation}
for all continua $E$ and $F$ in $D$ intersecting $\partial U$ and
$\partial V$. Here $M$ is the modulus, see (\ref{eq8.2.7}), and
$\Delta(E,F,D)$ the family of all paths
$\gamma:[a,b]\to\overline{{\Bbb R}^n}$ connecting $E$ and $F$ in
$D$, i.e., $\gamma(a)\in E$, $\gamma(b)\in F$ and $\gamma(t)\in D$
for all $t\in(a,b)$. We say that the boundary $\partial D$ is {\bf
weakly flat} if it is weakly flat at every point in $\partial D$.

We also say that a {\bf point} $x_0\in\partial D$ is {\bf strongly
accessible} if, for every neighborhood $U$ of the point $x_0$, there
exist a compactum $E$ in $D$, a neighborhood $V\subset U$ of $x_0$
and a number $\delta>0$ such that
\begin{equation}\label{eq1.6KR}M(\Delta(E,F,D))\ \geqslant\ \delta\end{equation} for all
continua $F$ in $D$ intersecting $\partial U$ and $\partial V$. We
say that the {\bf boundary} $\partial D$ is {\bf strongly
accessible} if every point $x_0\in\partial D$ is strongly
accessible.

\medskip

It is easy to see that if a domain $D$ in ${\Bbb R}^n$ is weakly
flat at a point $x_0\in\partial D$, then the point $x_0$ is strongly
accessible from $D$. Moreover, it was proved by us that if a domain
$D$ in ${\Bbb R}^n$ is weakly flat at a point $x_0\in\partial D$,
then $D$ is locally connected at $x_0$, see, e.g., Lemma 5.1 in
\cite{KR$_1$} or Lemma 3.15 in \cite{MRSY}.

\medskip

The notions of strong accessibility and weak flatness at boundary
points of a domain in ${\Bbb R}^n$ defined in \cite{KR}, see also
\cite{KR$_1$}, are localizations and generali\-zations of the
corresponding notions introduced in
\cite{MRSY$_5$}--\cite{MRSY$_6$}, cf. with properties $P_1$ and
$P_2$ by V\"ais\"al\"a in \cite{Va} and also with the quasiconformal
accessibility and the quasiconformal flatness by N\"akki in
\cite{Na$_1$}. Many theorems on a homeomorphic extension to the
boundary of quasiconformal mappings and their generalizations are
valid under the condition of weak flatness of boundaries. The
condition of strong accessibility plays a similar role for a
continuous extensi\-on of the mappings to the boundary.

\medskip

In the mapping theory and in the theory of differential equations,
it is often applied the so-called Lipschitz domains whose boundaries
are locally quasiconformal. Recall first that a map $\varphi:X\to Y$
between metric spaces $X$ and $Y$ is said to be {\bf Lipschitz}
provided ${\rm dist}(\varphi(x_1),\varphi(x_2))\leqslant M\cdot{\rm
dist}(x_1,x_2)$ for some $M<\infty$ and for all $x_1$ and $x_2\in
X$. The map $\varphi$ is called {\bf bi-Lipschitz} if, in addition,
$M^*{\rm dist}(x_1,x_2)\leqslant{\rm
dist}(\varphi(x_1),\varphi(x_2))$ for some $M^*>0$ and for all $x_1$
and $x_2\in X.$ Later on, $X$ and $Y$ are subsets of ${\Bbb R}^n$
with the Euclidean distance.

\medskip

It is said that a bounded domain $D$ in ${\Bbb R}^n$ is {\bf
Lipschitz} if every point $x_0\in\partial D$ has a neighborhood $U$
that can be mapped by a bi-Lipschitz homeomorphism $\varphi$ onto
the unit ball ${\Bbb B}^n\subset{\Bbb R}^n$ in such a way that
$\varphi(\partial D\cap U)$ is the intersection of ${\Bbb B}^n$ with
the a coordinate hyperplane and $f(x_0)=0$, see, e.g., \cite{Oht}.
Note that a bi-Lipschitz homeomorphism is quasiconformal and the
modulus is a quasiinvariant under such mappings. Hence the Lipschitz
domains have weakly flat boundaries. In particular, smooth and
convex bounded domains are so.

\bigskip

\cc
\section{On BMO, VMO and FMO functions}\label{5}

The BMO space was introduced by John and Nirenberg in \cite{JN} and
soon became one of the main concepts in harmonic analysis, complex
analysis, partial differential equations and related areas, see,
e.g., \cite{HKM} and \cite{ReRy}.

\medskip

Let $D$ be a domain in ${\Bbb R}^n$, $n\geqslant1$. Recall that a
real valued function $\varphi\in L^1_{\rm loc}(D)$ is said to be of
{\bf bounded mean oscillation} in $D$, abbr. $\varphi\in{\rm
BMO}(D)$ or simply $\varphi\in{\rm BMO}$, if
\begin{equation}\label{eq1.11} \Vert\varphi\Vert_*=
\sup\limits_{B\subset D}\ \ \ \dashint\limits_B\vert\varphi(z)-
\varphi_B\vert\,dm(z)<\infty\end{equation} where the supremum is
taken over all balls $B$ in $D$ and
\begin{equation}\label{eq1.12} \varphi_B=\dashint\limits_B
\varphi(z)\,dm(z)=\frac{1}{|B|}\int\limits_B\varphi(z)\,
dm(z)\end{equation} is the mean value of the function $\varphi$ over
$B$. Note that $L^{\infty}(D)\subset{\rm BMO}(D)\subset L^p_{\rm
loc}(D)$ for all $1\leqslant p<\infty$, see, e.g., \cite{ReRy}.

\medskip

A function $\varphi$ in BMO is said to have {\bf vanishing mean
oscillation}, abbr. $\varphi\in{\rm VMO}$, if the supremum in
(\ref{eq1.11}) taken over all balls $B$ in $D$ with
$|B|<\varepsilon$ converges to $0$ as $\varepsilon\to0$. VMO has
been introduced by Sarason in \cite{Sarason}. There are a number of
papers devoted to the study of partial differential equations with
coefficients of the class VMO, see, e.g., \cite{CFL}, \cite{ISbord},
\cite{MRV$^*$}, \cite{Pal} and \cite{Ra}.

\medskip

Following \cite{IR}, we say that a function $\varphi:D\to{\Bbb R}$
has {\bf finite mean oscillation at a point} $z_0\in D$, write
$\varphi\in{\rm FMO}(x_0)$, if \begin{equation}\label{FMO_eq2.4}
\overline{\lim\limits_{\varepsilon\to 0}}\ \ \ \mathchoice
{{\setbox0=\hbox{$\displaystyle{\textstyle -}{\int}$}
\vcenter{\hbox{$\textstyle -$}}\kern-.5\wd0}}
{{\setbox0=\hbox{$\textstyle{\scriptstyle -}{\int}$}
\vcenter{\hbox{$\scriptstyle -$}}\kern-.5\wd0}}
{{\setbox0=\hbox{$\scriptstyle{\scriptscriptstyle -}{\int}$}
\vcenter{\hbox{$\scriptscriptstyle -$}}\kern-.5\wd0}}
{{\setbox0=\hbox{$\scriptscriptstyle{\scriptscriptstyle -}{\int}$}
\vcenter{\hbox{$\scriptscriptstyle -$}}\kern-.5\wd0}}
\!\int_{B(z_0,\varepsilon)}
|{\varphi}(z)-\tilde{{\varphi}}_{\varepsilon}(z_0)|\,dm(z)<\infty
\end{equation} where \begin{equation}\label{FMO_eq2.5}
\tilde{{\varphi}}_{\varepsilon}(z_0)=\mathchoice
{{\setbox0=\hbox{$\displaystyle{\textstyle -}{\int}$}
\vcenter{\hbox{$\textstyle -$}}\kern-.5\wd0}}
{{\setbox0=\hbox{$\textstyle{\scriptstyle -}{\int}$}
\vcenter{\hbox{$\scriptstyle -$}}\kern-.5\wd0}}
{{\setbox0=\hbox{$\scriptstyle{\scriptscriptstyle -}{\int}$}
\vcenter{\hbox{$\scriptscriptstyle -$}}\kern-.5\wd0}}
{{\setbox0=\hbox{$\scriptscriptstyle{\scriptscriptstyle -}{\int}$}
\vcenter{\hbox{$\scriptscriptstyle -$}}\kern-.5\wd0}}
\!\int_{B(z_0,\varepsilon)}\varphi(z)\,dm(z)\end{equation} is the
mean value of the function $\varphi(z)$ over the ball
$B(z_0,\varepsilon)$. Condition (\ref{FMO_eq2.4}) includes the
assumption that $\varphi$ is integrable in some neighborhood of the
point $z_0$. By the triangle inequality, we obtain the following
statement.

\medskip

\begin{proposition}\label{FMO_pr2.1} {\it If for some collection of numbers
$\varphi_{\varepsilon}\in{\Bbb R}$,
$\varepsilon\in(0,\varepsilon_0]$,
\begin{equation}\label{FMO_eq2.7}
\overline{\lim\limits_{\varepsilon\to 0}}\ \ \ \mathchoice
{{\setbox0=\hbox{$\displaystyle{\textstyle -}{\int}$}
\vcenter{\hbox{$\textstyle -$}}\kern-.5\wd0}}
{{\setbox0=\hbox{$\textstyle{\scriptstyle -}{\int}$}
\vcenter{\hbox{$\scriptstyle -$}}\kern-.5\wd0}}
{{\setbox0=\hbox{$\scriptstyle{\scriptscriptstyle -}{\int}$}
\vcenter{\hbox{$\scriptscriptstyle -$}}\kern-.5\wd0}}
{{\setbox0=\hbox{$\scriptscriptstyle{\scriptscriptstyle -}{\int}$}
\vcenter{\hbox{$\scriptscriptstyle -$}}\kern-.5\wd0}}
\!\int_{B(z_0,\varepsilon)}
|{\varphi}(z)-{\varphi}_{\varepsilon}|\,dm(z)<\infty\,,
\end{equation} then $\varphi$ is of finite mean oscillation at
$z_0$.}
\end{proposition}

\medskip

Choosing in Proposition \ref{FMO_pr2.1}
$\varphi_{\varepsilon}\equiv0$, $\varepsilon\in(0,\varepsilon_0]$,
we have the following.

\medskip

\begin{corollary}\label{FMO_cor2.1} {\it If for a point $z_0\in D$
\begin{equation}\label{FMO_eq2.8}
\overline{\lim\limits_{\varepsilon\to 0}}\ \ \mathchoice
{{\setbox0=\hbox{$\displaystyle{\textstyle -}{\int}$}
\vcenter{\hbox{$\textstyle -$}}\kern-.5\wd0}}
{{\setbox0=\hbox{$\textstyle{\scriptstyle -}{\int}$}
\vcenter{\hbox{$\scriptstyle -$}}\kern-.5\wd0}}
{{\setbox0=\hbox{$\scriptstyle{\scriptscriptstyle -}{\int}$}
\vcenter{\hbox{$\scriptscriptstyle -$}}\kern-.5\wd0}}
{{\setbox0=\hbox{$\scriptscriptstyle{\scriptscriptstyle -}{\int}$}
\vcenter{\hbox{$\scriptscriptstyle -$}}\kern-.5\wd0}}
\!\int_{B(z_0,\varepsilon)}|{\varphi}(z)|\,dm(z)<\infty\,,\end{equation}
then $\varphi$ has finite mean oscillation at $z_0$}.
\end{corollary}

\medskip

Recall that a point $z_0\in D$ is called a {Lebesgue point} of a
function $\varphi:D\to{\Bbb R}$ if $\varphi$ is integrable in a
neighborhood of $z_0$ and \begin{equation}\label{FMO_eq2.7a}
\lim\limits_{\varepsilon\to 0}\ \ \ \mathchoice
{{\setbox0=\hbox{$\displaystyle{\textstyle -}{\int}$}
\vcenter{\hbox{$\textstyle -$}}\kern-.5\wd0}}
{{\setbox0=\hbox{$\textstyle{\scriptstyle -}{\int}$}
\vcenter{\hbox{$\scriptstyle -$}}\kern-.5\wd0}}
{{\setbox0=\hbox{$\scriptstyle{\scriptscriptstyle -}{\int}$}
\vcenter{\hbox{$\scriptscriptstyle -$}}\kern-.5\wd0}}
{{\setbox0=\hbox{$\scriptscriptstyle{\scriptscriptstyle -}{\int}$}
\vcenter{\hbox{$\scriptscriptstyle -$}}\kern-.5\wd0}}
\!\int_{B(z_0,\varepsilon)}|{\varphi}(z)-{\varphi}(z_0)|\,dm(z)=0\,.
\end{equation} It is known that almost every point in $D$ is a Lebesgue point for
every function $\varphi\in L^1(D)$. Thus, we have the following
conclusion.

\medskip

\begin{corollary}\label{FMO_cor2.7b} {\it Every function $\varphi:D\to{\Bbb R}$
which is locally integrable, has a finite mean oscillation at almost
every point in $D$.} \end{corollary}

\medskip

\begin{remark}\label{FMO_rmk2.13a} {\rm Note that the
function $\varphi(z)=\log(1/|z|)$ belongs to BMO in the unit disk
$\Delta$, see, e.g., \cite{ReRy}, p. 5, and hence also to FMO.
However, $\tilde{\varphi}_{\varepsilon}(0)\to\infty$ as
$\varepsilon\to0$, showing that the condition (\ref{FMO_eq2.8}) is
only sufficient but not necessary for a function $\varphi$ to be of
finite mean oscillation at $z_0$.} \end{remark}

\medskip

Clearly that ${\rm BMO}\subset{\rm FMO}\subset L^1_{\rm loc}$ but
FMO is not a subset of $L^p_{\rm loc}$ for any $p>1$, see examples
in Section 11.2 of the monograph \cite{MRSY}, in comparison with
${\rm BMO}_{\rm loc}\subset L^p_{\rm loc}$ for all $p\in[1,\infty)$.
Thus, ${\rm FMO}$ is essentially wider than ${\rm BMO}_{\rm loc}$.
The following lemma is key, see Corollary 2.3 in \cite{IR} and
Corollary 6.3 in \cite{MRSY}.

\medskip

\begin{lemma}\label{lem5.1log} {\it Let $D$ be a domain in ${\Bbb
R}^n$, $n\geqslant2$, $x_0\in D$, and let $\varphi:D\to{\Bbb R}$ be
a non-negative function of the class ${\rm FMO}(x_0)$. Then
\begin{equation}\label{eq5.1log}\int\limits_{\varepsilon<|x-x_0|<\varepsilon_0}
\frac{\varphi(x)\,dm(x)}{\left(|x-x_0|\log{\frac{1}{|x-x_0|}}\right)^n}
=O\left(\log{\log{\frac{1}{\varepsilon}}}\right)\quad\text{as}\quad\varepsilon\to0\end{equation}
for some $\varepsilon_0\in(0,\delta_0)$ where
$\delta_0=\min(e^{-e},d_0)$, $d_0=\sup\limits_{x\in
D}|x-x_0|$.}\end{lemma}

\bigskip

\cc
\section{On lower $Q$-homeomorphisms}

Let $\omega$ be an open set in $\overline{{\Bbb R}^k}$,
$k=1,\ldots,n-1$. Recall that a (continuous) mapping
$S:\omega\to{\Bbb R}^n$ is called a $k$-dimensional surface $S$ in
${\Bbb R}^n$. The number of preimages
\begin{equation}\label{eq8.2.3} N(S,y)={\rm card}\,S^{-1}(y)={\rm
card}\,\{x\in\omega:S(x)=y\},\ y\in{\Bbb R}^n\end{equation} is said
to be a {\bf multiplicity function} of the surface $S$. It is known
that the multiplicity function is lower semicontinuous, i.e.,
$$N(S,y)\geqslant\liminf_{m\to\infty}\:N(S,y_m)$$ for every sequence $y_m\in{\Bbb R}^n$, $m=1,2,\ldots\,$, such
that $y_m\to y\in{\Bbb R}^n$ as $m\to\infty$, see, e.g., \cite{RR},
p. 160. Thus, the function $N(S,y)$ is Borel measurable and hence
measurable with respect to every Hausdorff measure $H^k$, see, e.g.,
\cite{Sa}, p. 52.

\medskip

Recall that a $k$-dimensional Hausdorff area in ${\Bbb R}^n$ (or
simply {\bf area}) associated with a surface $S:\omega\to{\Bbb R}^n$
is given by \begin{equation}\label{eq8.2.4}
{\mathcal{A}}_S(B)={\mathcal{A}}^{k}_S(B):=\int\limits_B
N(S,y)\,dH^{k}y
\end{equation} for every Borel set $B\subseteq{\Bbb R}^n$ and,
more generally, for an arbitrary set that is measurable with respect
to $H^k$ in ${\Bbb R}^n$, cf. 3.2.1 in \cite{Fe} and 9.2 in
\cite{MRSY}.

\medskip

If $\varrho:{\Bbb R}^n\to\overline{{\Bbb R}^+}$ is a Borel function,
then its {\bf integral over} $S$ is defined by the equality
\begin{equation}\label{eq8.2.5} \int\limits_S \varrho\,d{\mathcal{A}}:=
\int\limits_{{\Bbb R}^n}\varrho(y)\,N(S,y)\,dH^ky\,.\end{equation}
Given a family $\Gamma$ of $k$-dimensional surfaces $S$, a Borel
function $\varrho:{\Bbb R}^n\to[0,\infty]$ is called {\bf
admissible} for $\Gamma$, abbr. $\varrho\in{\rm adm}\,\Gamma$, if
\begin{equation}\label{eq8.2.6}\int\limits_S\varrho^k\,d{\mathcal{A}}\geqslant1\end{equation}
for every $S\in\Gamma$. The {\bf modulus} of $\Gamma$ is the
quantity
\begin{equation}\label{eq8.2.7} M(\Gamma)=\inf_{\varrho\in{\rm adm}\,\Gamma}
\int\limits_{{\Bbb R}^n}\varrho^n(x)\,dm(x)\,.\end{equation} We also
say that a Lebesgue measurable function $\varrho:{\Bbb
R}^n\to[0,\infty]$ is {\bf extensively admissible} for a family
$\Gamma$ of $k$-dimensional surfaces $S$ in ${\Bbb R}^n$, abbr.
$\varrho\in{\rm ext\, adm}\,\Gamma$,  if a subfamily of all surfaces
$S$ in $\Gamma$, for which (\ref{eq8.2.6}) fails, has the modulus
zero.

Given domains $D$ and $D'$ in $\overline{{\Bbb R}^n}={\Bbb
R}^n\cup\{\infty\}$, $n\geqslant2$,
$x_0\in\overline{D}\setminus\{\infty\}$, and a measurable function
$Q:{{\Bbb R}^n}\to(0,\infty)$, we say that a homeomorphism $f:D\to
D'$ is a {\bf lower $Q$-homeomorphism at the point} $x_0$ if
\begin{equation}\label{eqOS1.10} M(f\Sigma_{\varepsilon})\geqslant
\inf\limits_{\varrho\in{\rm
ext\,adm}\,\Sigma_{\varepsilon}}\int\limits_{D\cap
R_{\varepsilon}}\frac{\varrho^n(x)}{Q(x)}\,dm(x)\end{equation} for
every ring $R_{\varepsilon}=\{x\in{\Bbb
R}^n:\varepsilon<|x-x_0|<\varepsilon_0\}\,,\quad\varepsilon\in(0,\varepsilon_0)\,,\
\varepsilon_0\in(0,d_0)$, where $d_0=\sup\limits_{x\in D}|x-x_0|$,
and $\Sigma_{\varepsilon}$ denotes the family of all intersections
of the spheres $S(x_0,r)=\{x\in{\Bbb R}^n:|x-x_0|=r\}\,,
r\in(\varepsilon,\varepsilon_0)\,,$ with $D$.

\medskip

The notion of lower $Q$-homeomorphisms in the standard way can be
extended to $\infty$ through inversions. We also say that a
homeomorphism $f:D\to{\overline{{\Bbb R}^n}}$ is a {\bf lower
$Q$-homeomor\-phism on} $\partial D$ if $f$ is a lower
$Q$-homeomorphism at every point $x_0\in\partial{D}$.

\medskip

We proved the following significant statements on lower
$Q$-ho\-meo\-mor\-phisms, see Theorem 10.1 (Lemma 6.1) in
\cite{KR$_1$} or Theorem 9.8 (Lemma 9.4) in \cite{MRSY}.

\bigskip

\begin{proposition}\label{prKR2.1} {\it Let $D$ and $D'$ be bounded domains in
${\Bbb R}^n$, $n\geqslant2$, $Q:D\to(0,\infty)$ a measurable
function and $f:D\to D'$ a lower $Q$-homeomorphism on $\partial D$.
Suppose that the domain $D$ is locally connected on $\partial D$ and
that the domain $D'$ has a (strongly accessible) weakly flat
boundary. If
\begin{equation}\label{eqKPR2.11}\int\limits_{0}^{\delta(x_0)} \frac{dr}{||\,Q||_{n-1}(x_0,r)}\
=\ \infty\qquad\forall\ x_0\in\partial D\end{equation} for some
$\delta(x_0)\in(0,d(x_0))$ where $d(x_0)=\sup\limits_{x\in
D}\,|\,x-x_0|$ and
$$||Q||_{n-1}(x_0,r)=\left(\int\limits_{D\cap
S(x_0,r)}Q^{n-1}(x)\,d{\cal A}\right)^{\frac{1}{n-1}}\,,$$ then $f$
can be extended to a (continuous) homeomorphic mapping
$\overline{f}:\overline{D}\to\overline{D'}$.}\end{proposition}

\medskip

\section{A connection with the Orlicz--Sobolev classes}

Given a mapping $f:D \to \Rn$ with partial derivatives a.e., recall
that $f^\prime(x)$ denotes the Jacobian matrix of $f$ at $x \in D$
if it exists, $J(x)=J(x,f)=\det f^\prime(x)$ is the Jacobian of $f$
at $x$, and $\Vert f^\prime(x)\Vert$ is the operator norm of
$f^\prime(x)$, i.e.,
\begin{equation} \label{eq4.1.2}
\Vert f^\prime(x)\Vert =\max \{|f^\prime(x)h|: h \in \Rn, |h|=1\}.
\end{equation}
We also let
\begin{equation} \label{eq4.1.3}
l(f^\prime(x))= \min \{ |f^\prime(x)h|: h \in \Rn, |h|=1\}.
\end{equation}
The {\bf outer dilatation} of $f$ at $x$ is defined by
\begin{equation} \label{eq4.1.4} K_O(x)=K_O(x,f)= \left
\{\begin{array}{rl}
\frac{\Vert f^\prime(x)\Vert^n}{|J(x,f)|} & {\rm if } \ J(x,f) \neq 0, \\
1 & {\rm if} \ f^\prime(x)=0, \\ \infty & {\rm } \text{otherwise},
\end{array} \right. \end{equation} the {\bf inner dilatation} of $f$ at $x$
by \begin{equation} \label{eq4.1.5} K_I(x)=K_I(x,f)= \left
\{\begin{array}{rl}
\frac{|J(x,f)|}{l(f^\prime(x))^n} &{\rm if } \ J(x,f) \neq 0, \\
1 & {\rm if } \ f^\prime(x)=0, \\
\infty &{\rm} \text{otherwise},
\end{array} \right. \end{equation} Note that, see, e.g.,
Section 1.2.1 in \cite{Re},
\begin{equation} \label{eq4.1.444} K_O(x,f)\leq K^{n-1}_I(x,f)\ \ \ \ \
\mbox{and}\ \ \ \ \ K_I(x,f) \leq K^{n-1}_O(x,f)\ ,\end{equation} in
particular, $K_O(x,f) < \infty $ a.e. if and only if $K_I(x,f) <
\infty $  a.e. The latter is equivalent to the condition that a.e.
either $ \det f^\prime(x)
> 0$ or $f^\prime(x)=0.$

\medskip

Now, recall that a homeomorphism $f$ between domains $D$ and $D'$ in
${\Bbb R}^n$, $n\geqslant2$, is called of {\bf finite distortion} if
$f\in W^{1,1}_{\rm loc}$ and \begin{equation}\label{eqOS1.3} \Vert
f'(x)\Vert^n\leqslant K(x)\cdot J_f(x)\end{equation} with some a.e.
finite function $K$. The term is due to Tadeusz Iwaniec. In other
words, (\ref{eqOS1.3}) just means that dilatations $K_O(x,f)$ and
$K_I(x,f)$ are finite a.e.

\medskip

In view of (\ref{eq4.1.444}), the next statement says on a stronger
modulus estimate than the obtained in \cite{KRSS1}, Theorem 4.1, in
terms of the outer dilatation $K_O(x,f)$. It is key for deriving
consequences from our theory of lower $Q-$ho\-meo\-mor\-phisms.

\medskip

\begin{theorem}\label{thOS4.1} {\it
Let $D$ and $D'$ be domains in ${\Bbb R}^n$, $n\geqslant3$, and let
$\varphi:{\Bbb R}^+\to{\Bbb R}^+$ be a nondecreasing function such
that, for some $t_*\in{\Bbb R}^+$,
\begin{equation}\label{eqOS4.1}
\int\limits_{t_*}^{\infty}\left[\frac{t}{\varphi(t)}\right]^
{\frac{1}{n-2}}dt<\infty\,.\end{equation} Then each homeomorphism
$f:D\to D'$ of finite distortion in the class $W^{1,\varphi}_{\rm
loc}$ is a lower $Q$-homeomorphism at every point
$x_0\in\overline{D}$ with $Q(x)=\left[
K_I(x,f)\right]^{\frac{1}{n-1}}$.}
\end{theorem}

\medskip

\begin{proof}
Let $B$ be a (Borel) set of all points $x\in D$ where $f$ has a
total differential $f'(x)$ and $J_f(x)\ne0$. Then, applying
Kirszbraun's theorem and uni\-que\-ness of approximate differential,
see, e.g., 2.10.43 and 3.1.2 in \cite{Fe}, we see that $B$ is the
union of a countable collection of Borel sets $B_l$,
$l=1,2,\ldots\,$, such that $f_l=f|_{B_l}$ are bi-Lipschitz
homeomorphisms, see, e.g., 3.2.2 as well as 3.1.4 and 3.1.8 in
\cite{Fe}. With no loss of generality, we may assume that the $B_l$
are mutually disjoint. Denote also by $B_*$ the rest of all points
$x\in D$ where $f$ has the total differential but with $f'(x)=0$.

By the construction the set $B_0:=D\setminus\left(B\bigcup
B_*\right)$ has Lebesgue measure zero, see Theorem 1 in \cite{KRSS}.
Hence ${\mathcal A}_S(B_0)=0$ for a.e. hypersurface $S$ in ${\Bbb
R}^n$  and, in particular, for a.e. sphere $S_r:=S(x_0,r)$ centered
at a prescribed point $x_0\in\overline{D}$, see Theorem 2.11 in
\cite{KR$_7$} or Theorem 9.1 in \cite{MRSY}. Thus, by Corollary 4 in
\cite{KRSS} ${\mathcal A}_{S_r^*}(f(B_0))=0$ as well as ${\mathcal
A}_{S_r^*}(f(B_*))=0$ for a.e. $S_r$ where $S_r^*=f(S_r)$.

Let $\Gamma$ be the family of all intersections of the spheres
$S_r$, $r\in(\varepsilon,\varepsilon_0)$,
$\varepsilon_0<d_0=\sup\limits_{x\in D}\,|x-x_0|$, with the domain
$D$. Given $\varrho_*\in{\rm adm}\,f(\Gamma)$ such that
$\varrho_*\equiv0$ outside of $f(D)$, set $\varrho\equiv0$ outside
of $D$ and on $D\setminus B$ and, moreover,
$$
\varrho(x):=\Lambda(x)\cdot\varrho_*(f(x))\qquad{\rm for}\ x\in B
$$
where
$$
\Lambda(x)\ :=\ \left[\ J_f(x)\cdot K^{\frac{1}{n-1}}_I(x,f)\
\right]^{\frac{1}{n}}\ =\ \left[\ \frac{\mbox{det}\,
f^{\prime}(x)}{l( f^{\prime}(x))}\ \right]^{\frac{1}{n-1}}\ =$$
$$=\ \left[\
\lambda_2\cdot\ldots\cdot\lambda_n\ \right]^{\frac{1}{n-1}}\
\geqslant \ \left[\ J_{n-1}(x)\ \right]^{\frac{1}{n-1}}\qquad{\rm
for\ a.e.}\ x\in B\ ;
$$
here as usual $\lambda_n\geqslant\ldots\geqslant\lambda_1$ are
principal dilatation coefficients of $f^{\prime}(x)$, see, e.g.,
Section I.4.1 in \cite{Re}, and $J_{n-1}(x)$ is the
$(n-1)-$dimensional Jacobian of $f|_{S_r}$ at $x$  where
$r=|x-x_0|$, see Section 3.2.1 in \cite{Fe}.

Arguing piecewise on $B_l$, $l=1,2,\ldots\,$, and taking into
account Kirszbraun's theorem, by Theorem 3.2.5 on the change of
variables  in \cite{Fe}, we have that
$$\int\limits_{S_r}\varrho^{n-1}\,d{\mathcal A}\geqslant\int\limits_{S_{*}^r}\varrho_{*}^{n-1}\,d{\mathcal A}\geqslant1$$
for a.e. $S_r$ and, thus, $\varrho\in{\rm ext\,adm}\,\Gamma$.

The change of variables on each $B_l$, $l=1,2,\ldots\,$, see again
Theorem 3.2.5 in \cite{Fe}, and countable additivity of integrals
give also the estimate
$$\int\limits_{D}\frac{\varrho^n(x)}{K^{\frac{1}{n-1}}_I(x)}\,dm(x)\leqslant
\int\limits_{f(D)}\varrho^n_*(x)\,dm(x)$$ and the proof is complete.
\end{proof}$\Box$

\medskip

\begin{corollary}\label{corOS4.1} {\it
Each homeomorphism $f$ with finite distortion in ${\Bbb R}^n$,
$n\geqslant3$, of the class $W^{1,p}_{\rm loc}$ for $p>n-1$ is a
lower $Q$-homeomorphism at every point $x_0\in\overline{D}$ with
$Q=K^{\frac{1}{n-1}}_I$.}
\end{corollary}

\medskip

\section{Boundary behavior of Orlicz--Sobolev classes}

In this section  we assume that $\varphi:{\Bbb R}^+\to{\Bbb R}^+$ is
a nondecreasing function such that, for some $t_*\in{\Bbb R}^+$,
\begin{equation}\label{eqOSKRSS}\int\limits_{t_*}^{\infty}\left[\frac{t}{\varphi(t)}\right]^
{\frac{1}{n-2}}dt<\infty\,.\end{equation} The continuous extension
to the boundary of the inverse mappings has a simpler criterion than
for the direct mappings. Hence we start from the former. Namely, in
view of Theorem \ref{thOS4.1}, we have by Theorem 9.1 in
\cite{KR$_1$} or Theorem 9.6 in \cite{MRSY} the next statement.

\medskip

\begin{theorem}\label{thKPR8.2} {\it Let $D$ and $D'$ be bounded domains in ${\Bbb
R}^n$, $n\geqslant3$, $D$ be locally connected on $\partial D$ and
$\partial D'$ be weakly flat. Suppose that $f$ is a homeo\-morphism
of $D$ onto $D'$ in a class $W^{1,\varphi}_{\rm loc}$ with condition
(\ref{eqOSKRSS}) and $K_I\in L^{1}(D)$. Then $f^{-1}$ can be
extended to a continuous mapping of $\overline{D^{\prime}}$ onto
$\overline{D}$.}
\end{theorem}

\medskip

However, as it follows from the example in Proposition 6.3 in
\cite{MRSY}, see also (\ref{eq4.1.444}), any degree of integrability
$K_I\in L^q(D)$, $q\in [1,\infty)$, cannot guarantee  the extension
by continuity to the boundary of the direct mappings.

\medskip

Also by Theorem \ref{thOS4.1}, we have the following consequence of
Proposition \ref{prKR2.1}.

\medskip

\begin{theorem}\label{thKR9.111} {\it Let $D$ and $D'$ be bounded domains in
${\Bbb R}^n$, $n\geqslant3$, $D$ be locally connected on $\partial
D$ and $\partial D'$ be (strongly accessible) weakly flat. Suppose
that $f:D\to D'$ is a homeomorphism of finite distortion in
$W^{1,\varphi}_{\rm loc}$ with condition (\ref{eqOSKRSS}) such that
\begin{equation}\label{eqKPR2.1}
\int\limits_{0}^{\delta(x_0)}\frac{dr}{||K_I||^{\frac{1}{n-1}}(x_0,r)}
=\infty\qquad\forall\ x_0\in\partial D\end{equation} for some
$\delta(x_0)\in(0,d(x_0))$ where $d(x_0)=\sup\limits_{x\in
D}|x-x_0|$ and
$$||K_I||(x_0,r)=\int\limits_{D\cap S(x_0,r)}K_I(x,f)\ d{\mathcal A}\ .$$
Then  $f$ can be extended to a (continuous) homeomorphic map
$\overline{f}:\overline{D}\to\overline{D'}$.} \end{theorem}

\medskip

In particular, as a consequence of Theorem \ref{thKR9.111}, we
obtain the following genera\-lization of the well-known theorems of
Gehring--Martio and Martio--Vuorinen on a homeomorphic extension to
the boundary of quasiconformal mappings between the so--called QED
domains, see \cite{GM} and \cite{MV}.

\medskip

\begin{corollary}\label{thKPR9.2} {\it Let $D$ and $D'$ be bounded domains in
${\Bbb R}^n$, $n\geqslant3$, $D$ be locally connected on $\partial
D$ and $\partial D'$ be (strongly accessible) weakly flat. Suppose
that $f:D\to D'$ is a homeo\-morphism of finite distortion in the
class $W^{1,p}_{\rm loc}$, $p>n-1$. If (\ref{eqKPR2.1}) holds, then
$f$ can be extended to a (continuous) homeomorphic map
$\overline{f}:\overline{D}\to\overline{D'}$.}
\end{corollary}

\medskip

\begin{lemma}\label{lemOSKRSS12.1} {\it Let $D$ and $D'$ be bounded domains in
${\Bbb R}^n$, $n\geqslant3$, $D$ be locally connected on $\partial
D$ and $\partial D'$ be (strongly accessible) weakly flat. Suppose
that $f:D\to D'$ is a homeomorphism of finite distortion in
$W^{1,\varphi}_{\rm loc}$ with condition (\ref{eqOSKRSS}) such that
\begin{equation}\label{omal} \int\limits_{D(x_0,\varepsilon,\varepsilon_0)}
K_I(x,f)\cdot\psi^n_{x_0,\varepsilon}(|x-x_0|)\,dm(x)=o(I_{x_0}^n(\varepsilon))\
{\rm as}\ \varepsilon\to0\ \forall\ x_0\in\partial D\end{equation}
where $D(x_0,\varepsilon,\varepsilon_0)=\{x\in
D:\varepsilon<|x-x_0|<\varepsilon_0\}$ for some
$\varepsilon_0\in(0,\delta_0)$, $\delta_0=\delta(x_0)=\sup_{x\in
D}|x-x_0|$, and $\psi_{x_0,\varepsilon}(t)$ is a family of
non-negative measurable (by Lebesgue) functions on $(0,\infty)$ such
that \begin{equation}\label{eq5.3} 0\ <\ I_{x_0}(\varepsilon)\ =\
\int\limits_{\varepsilon}^{\varepsilon_0} \psi_{x_0,\varepsilon}(t)\
dt\ <\ \infty\qquad\forall\ \varepsilon\in(0,\varepsilon_0)\
.\end{equation} Then  $f$ can be extended to a (continuous)
homeomorphic map $\overline{f}:\overline{D}\to\overline{D'}$. }
\end{lemma}

\medskip

\begin{proof} Lemma \ref{lemOSKRSS12.1} is an immediate consequence of Theorems
\ref{thOS4.1} and \ref{thKR9.111} ta\-king into account Lemma 3.7 in
the work \cite{RS}, see Lemma 7.4 in the monograph \cite{MRSY}, and
extending $K_I(x,f)$ by zero outside of $D.$ \end{proof}$\Box$

\medskip

Choosing in Lemma \ref{lemOSKRSS12.1} $\psi(t)=1/(t\log{1/t})$ and
applying Lemma \ref{lem5.1log}, we obtain the following result.

\medskip

\begin{theorem}\label{thKPRS12a*} {\it Let $D$ and $D'$ be bounded domains in
${\Bbb R}^n$, $n\geqslant3$, $D$ be locally connected on $\partial
D$ and $\partial D'$ be (strongly accessible) weakly flat. Suppose
that $f:D\to D'$ is a homeomorphism in $W^{1,\varphi}_{\rm loc}$
with condition (\ref{eqOSKRSS}) such that
\begin{equation}\label{eqKPRS12a*}K_I(x,f)\ \leqslant\ Q(x)\quad\quad\quad{\rm
a.e.\ in}\ D\end{equation} for a function $Q:{\Bbb R}^n\to{\Bbb
R}^n$, $Q\in{\rm FMO}(x_0)$ for all $x_0\in\partial D$. Then  $f$
can be extended to a (continuous) homeomorphic map
$\overline{f}:\overline{D}\to\overline{D'}$.}
\end{theorem}

\medskip

In the following consequences, we assume that $K_I(x,f)$ is extended
by zero outside of $D$.

\medskip

\begin{corollary}\label{corKPRS12a*} {\it
 In particular, the conclusions
of Theorem \ref{thKPRS12a*} hold if
\begin{equation}\label{eqKPRS12b*}\overline{\lim\limits_{\varepsilon\to0}}\
\ \dashint_{B(x_0,\varepsilon)}K_I(x,f)\ dm(x)\ <\
\infty\qquad\forall\ x_0\in\partial D\ .\end{equation}}
\end{corollary}

\medskip

Similarly, choosing in Lemma \ref{lemOSKRSS12.1} the function
$\psi(t)=1/t$, we come to the following more general statement.

\medskip

\begin{theorem}\label{thKPRS12b*} {\it Let $D$ and $D'$ be bounded domains in
${\Bbb R}^n$, $n\geqslant3$, $D$ be locally connected on $\partial
D$ and $\partial D'$ be (strongly accessible) weakly flat. Suppose
that $f:D\to D'$ is a homeomorphism in $W^{1,\varphi}_{\rm loc}$
with condition (\ref{eqOSKRSS}) such that
\begin{equation}\label{eqKPRS12c*}
\int\limits_{\varepsilon<|x-x_0|<\varepsilon_0}K_I(x,f)\
\frac{dm(x)}{|x-x_0|^n}\ =\
o\left(\left[\log\frac{\varepsilon_0}{\varepsilon}\right]^n\right)\qquad\forall\
x_0\in\partial D\end{equation} as $\varepsilon\to0$ for some
$\varepsilon_0\in(0,\delta_0)$ where
$\delta_0=\delta(x_0)=\sup_{x\in D}|x-x_0|$. Then  $f$ can be
extended to a (continuous) homeomorphic map
$\overline{f}:\overline{D}\to\overline{D'}$.}
\end{theorem}

\medskip

\begin{corollary}\label{corKPRS12b*} {\it
 The condition (\ref{eqKPRS12c*}) and the
assertion of Theorem \ref{thKPRS12b*} hold if
\begin{equation}\label{eqKPRS12d*}K_I(x,f)\ =\ o\left(\left[\log\frac{1}{|x-x_0|}\right]^{n-1}\right)\end{equation}
as $x\to x_0$. The same holds if \begin{equation}\label{eqKPRS12e*}
k_f(r)=o\left(\left[\log\frac{1}{r}\right]^{n-1}\right)
\end{equation} as $r\to0$ where $k_f(r)$ is the mean value of the
function $K_I(x,f)$ over the sphere $|x-x_0|=r$.}
\end{corollary}

\medskip

\begin{remark}\label{rmKRRSa*} Choosing in Lemma \ref{lemOSKRSS12.1} the function
$\psi(t)=1/(t\log{1/t})$ instead of $\psi(t)=1/t$, we are able to
replace (\ref{eqKPRS12c*}) by
\begin{equation}\label{eqKPRS12f*}
\int\limits_{\varepsilon<|x-x_0|<1}\frac{K_I(x,f)\
dm(x)}{\left(|x-x_0|\log{\frac{1}{|x-x_0|}}\right)^n}
=o\left(\left[\log\log\frac{1}{\varepsilon}\right]^n\right)\end{equation}
and (\ref{eqKPRS12e*}) by
\begin{equation}\label{eqKPRS12g*}
k_f(r)=o\left(\left[\log\frac{1}{r}\log\log\frac{1}{r}\right]^{n-1}\right).\end{equation}
Thus, it is sufficient to require that
\begin{equation}\label{eqKPRS12h*}k_f(r)=O\left(\left[\log\frac{1}{r}\right]^{n-1}\right)
\end{equation}

In general, we could give here the whole scale of the corresponding
conditions in terms of $\log$ using functions $\psi(t)$ of the form
$1/(t\,\log\ldots\log1/t)$. \end{remark}

\medskip

\begin{theorem}\label{thKR4.1} {\it Let $D$ and $D'$ be bounded domains in
${\Bbb R}^n$, $n\geqslant3$, $D$ be locally connected on $\partial
D$ and $\partial D'$ be (strongly accessible) weakly flat. Suppose
that $f:D\to D'$ is a homeomorphism in $W^{1,\varphi}_{\rm loc}$
with condition (\ref{eqOSKRSS}) such that
\begin{equation}\label{eqKR4.1}
\int\limits_{D}\Phi(K_I(x,f))\ dm(x)\ <\ \infty\end{equation} for a
non-decreasing convex function $\Phi:\overline{{\Bbb
R}^+}\to\overline{{\Bbb R}^+}$. If, for some $\delta>\Phi(0)$,
\begin{equation}\label{eqKR4.2}\int\limits_{\delta}^{\infty}\frac{d\tau}{\tau\left[\Phi^{-1}(\tau)\right]^{\frac{1}{n-1}}}=
\infty\end{equation} then  $f$ can be extended to a (continuous)
homeomorphic map $\overline{f}:\overline{D}\to\overline{D'}$.}
\end{theorem}

\medskip

Indeed, by Theorem 3.1 and Corollary 3.2 in
 \cite{RSY}, (\ref{eqKR4.1}) and
(\ref{eqKR4.2}) imply (\ref{eqKPR2.1}) and, thus, Theorem
\ref{thKR4.1} is a direct consequence of Theorem \ref{thKR9.111}.

\medskip

\begin{corollary}\label{corOSKRSS4.6.3} {\it The conclusion of Theorem
\ref{thKR4.1} holds if, for some $\alpha>0$,
\begin{equation}\label{eqOSKRSS4.6.6}
\int\limits_{D}e^{\alpha K_I(x,f)}\ dm(x)\ <\ \infty\
.\end{equation} }
\end{corollary}

\medskip

\begin{remark}\label{rmKR4.1} Note that by Theorem 5.1 and Remark 5.1 in
\cite{KR$_3$} the conditions (\ref{eqKR4.2}) are not only sufficient
but also necessary for continuous extension to the boundary of $f$
with the integral constraint (\ref{eqKR4.1}).

\medskip

Recall that by Theorem 2.1 in \cite{RSY}, see also Proposition 2.3
in \cite{RS1}, condition (\ref{eqKR4.2}) is equivalent to a series
of other conditions and, in particular, to the following condition
\begin{equation}\label{eqKR4.4}
\int\limits_{\delta}^{\infty}\log\,\Phi(t)\,\frac{dt}{t^{n'}}=+\infty\end{equation}
for some $\delta>0$ where $\frac{1}{n'}+\frac{1}{n}=1$, i.e., $n'=2$
for $n=2$, $n'$ is strictly decreasing in $n$ and $n'=n/(n-1)\to1$
as $n\to\infty$. \end{remark}

\medskip

Finally, note that all these results hold, for instance, if $f\in
W^{1,p}_{\rm loc}$, $p>n-1$.

\medskip

\section{On finitely bi--Lipschitz mappings}

Given an open set $\Omega\subseteq \Bbb{R}^n$, $n\geqslant 2$,
following Section 5 in \cite{KR$_7$}, see also Section 10.6 in
\cite{MRSY}, we say that a mapping $f:\Omega\to \Bbb{R}^n$ is {\bf
finitely bi-Lipschitz} if
\begin{equation}\label{eq8.12.2} 0\ <\ l(x,f)\ \leqslant\ L(x,f)\ <\
\infty \ \ \ \ \ \forall\ x\in\Omega\end{equation} where
\begin{equation} \label{eq8.1.6} L(x,f)\ =\
\limsup_{y\to x}\ \frac{|f(y)-f(x)|}{|y-x|}
\end{equation} and
\begin{equation}\label{eq8.1.7} l(x,f)\ =\ \liminf_{y\to
x}\ \frac{|f(y)-f(x)|}{|y-x|}\ .\end{equation}

\medskip

By the classic Rademacher--Stepanov theorem, we obtain from the
right hand inequality in (\ref{eq8.12.2}) that finitely bi-Lipschitz
mappings are differentiable a.e. and from the left hand inequality
in (\ref{eq8.12.2}) that $J_f(x)\ne 0$ a.e. Moreover, such mappings
have $(N)-$property with respect to each Hausdorff measure, see,
e.g., either Lemma 5.3 in \cite{KR$_7$} or Lemma 10.6 \cite{MRSY}.
Thus, the proof of the following theorems is perfectly similar to
one of Theorem \ref{thOS4.1} and hence we omit it, cf. also similar
but weaker Corollary 5.15 in \cite{KR$_7$} and Corollary 10.10 in
\cite{MRSY} formulated in terms of the outer dilatation $K_O$.

\medskip

\begin{theorem}\label{pr8.12.15} {\it Every finitely bi-Lipschitz
homeomorphism $f:\Omega\to \Bbb{R}^n$, $n\geqslant 2$, is a lower
$Q$-homeomorphism with $Q=K^{\frac{1}{n-1}}_I$.}
\end{theorem}

\medskip

All results for finitely bi-Lipschitz homeomorphisms are perfectly
si\-mi\-lar to the corresponding results for homeomorphisms with
finite distortion in the Orlich--Sobolev classes from the last
section. Hence we will not formulate all them in the explicit form
here in terms of inner dilatation $K_I$.

\medskip

We give here for instance only one of these results.

\medskip

\begin{theorem}\label{thKR4.1fb} {\it Let $D$ and $D'$ be bounded domains in
${\Bbb R}^n$, $n\geqslant2$, $D$ be locally connected on $\partial
D$ and $\partial D'$ be (strongly accessible) weakly flat. Suppose
that $f:D\to D'$ is a finitely bi-Lipschitz homeomorphism such that
\begin{equation}\label{eqKR4.1fb}
\int\limits_{D}\Phi(K_I(x,f))\ dm(x)\ <\ \infty\end{equation} for a
non-decreasing convex function $\Phi:\overline{{\Bbb
R}^+}\to\overline{{\Bbb R}^+}$. If, for some $\delta>\Phi(0)$,
\begin{equation}\label{eqKR4.2fb}\int\limits_{\delta}^{\infty}\frac{d\tau}{\tau\left[\Phi^{-1}(\tau)\right]^{\frac{1}{n-1}}}=
\infty\end{equation} then $f$ can be extended to a (continuous)
homeomorphic map $\overline{f}:\overline{D}\to\overline{D'}$.}
\end{theorem}

\medskip

\begin{corollary}\label{corOSKRSS4.6.3fb} {\it The conclusion of Theorem
\ref{thKR4.1fb} holds if, for some $\alpha>0$,
\begin{equation}\label{eqOSKRSS4.6.6fb}
\int\limits_{D}e^{\alpha K_I(x,f)}\ dm(x)\ <\ \infty\
.\end{equation} }
\end{corollary}

\medskip
\noindent
{\bf Denis Kovtonyuk and Vladimir Ryazanov,}\\
Institute of Applied Mathematics and Mechanics,\\
National Academy of Sciences of Ukraine,\\
74 Roze Luxemburg Str., Donetsk, 83114, Ukraine,\\
denis$\underline{\ \ }$\,kovtonyuk@bk.ru, vl.ryazanov1@gmail.com

\end{document}